\documentclass[12pt,twoside]{article}

\usepackage{amssymb}
\usepackage{amsmath}
\usepackage{bbm}
\usepackage{mathrsfs}
\usepackage{float}
\usepackage{graphicx}
\usepackage{multicol}
\usepackage{pifont}

\newcommand{\br}{ }
\newcommand{\brr}{, }

\makeatletter
\gdef\th@mychange{\normalfont\slshape
   \def\@begintheorem##1##2{\item
        [\hskip\labelsep \theorem@headerfont ##2. ##1  \,--\!--\!--\!--  ]}%
 \def\@opargbegintheorem##1##2##3{%
   \item[\hskip\labelsep \theorem@headerfont ##2. ##1\ {\upshape(}##3{\upshape)}. \,-----  ]}}
\makeatother

\RequirePackage{theorem}
\theoremstyle{mychange}

{\theorembodyfont{\rmfamily}\newtheorem{ttt}{}[section]}
{\theorembodyfont{\rmfamily}\newtheorem{defi}[ttt]{Definition.}}
{\theorembodyfont{\rmfamily}\newtheorem{rem}[ttt]{Remark.}}
{\theorembodyfont{\rmfamily}\newtheorem{rems}[ttt]{Remarks.}}
{\theorembodyfont{\rmfamily}\newtheorem{exa}[ttt]{Example.}}
{\theorembodyfont{\rmfamily}\newtheorem{exas}[ttt]{Examples.}}
{\theorembodyfont{\rmfamily}\newtheorem{conv}[ttt]{Convention.}}
{\theorembodyfont{\rmfamily}\newtheorem{qu}[ttt]{Question.}}
{\theorembodyfont{\rmfamily}\newtheorem{qus}[ttt]{Questions.}}
{\theorembodyfont{\rmfamily}}

{\theorembodyfont{\rmfamily}}
{\theorembodyfont{\rmfamily}\newtheorem{exo}[ttt]{Example}}

{\theorembodyfont{\itshape}\newtheorem{fac}[ttt]{Fact.}}
{\theorembodyfont{\itshape}\newtheorem{lem}[ttt]{Lemma.}}
{\theorembodyfont{\itshape}\newtheorem{prop}[ttt]{Proposition.}}
{\theorembodyfont{\itshape}\newtheorem{coro}[ttt]{Corollary.}}

{\theorembodyfont{\itshape}}
{\theorembodyfont{\itshape}\newtheorem{lemo}[ttt]{Lemma}}
{\theorembodyfont{\itshape}\newtheorem{propo}[ttt]{Proposition}}
{\theorembodyfont{\itshape}\newtheorem{cono}[ttt]{Conjecture}}

{\theorembodyfont{\rmfamily}\newtheorem{remo}[ttt]{Remark}}

{\theorembodyfont{\rmfamily}\newtheorem{lui}[ttt]{The method of van Luijk in detail.}}

\renewcommand{\atop}[2]{\genfrac{}{}{0pt}{}{#1}{#2}}

\newcommand{\calO}{\mathscr{O}}

\newcommand{\calV}{\mathscr{V}}

\newcommand{\bbC}{{\mathbbm C}}
\newcommand{\bbF}{{\mathbbm F}}
\newcommand{\bbG}{{\mathbbm G}}

\newcommand{\bbQ}{{\mathbbm Q}}
\newcommand{\bbZ}{{\mathbbm Z}}

\newcommand{\id}{\mathop{\rm id}}
\newcommand{\Tr}{\mathop{\rm Tr}}
\newcommand{\diag}{\mathop{\rm diag}}
\newcommand{\bP}{{\bf P}}
\newcommand{\rk}{\mathop{\text{\rm rk}}}
\newcommand{\Gal}{\mathop{\rm Gal}}
\newcommand{\Frob}{\mathop{\rm Frob}}
\newcommand{\Pic}{\mathop{\rm Pic}}
\newcommand{\NS}{\mathop{\rm NS}}
\newcommand{\End}{\mathop{\rm End}}
\newcommand{\Br}{\mathop{\rm Br}}

\newcommand{\reg}{\text{reg}}
\newcommand{\alg}{\text{alg}}
\renewcommand{\top}{\text{top}}
\newcommand{\et}{\text{\rm \'et}}

\newcommand{\bigtimes}{\mathop{\mbox{\Pisymbol{pzd}{53}}}\limits}

\newcounter{abc}
\newenvironment{abc}{\begin{list}{\rm \alph{abc}) }%
{\usecounter{abc} \leftmargin=0.0pt \labelsep=0.0pt %
\listparindent=0.0pt \labelwidth=0.0pt \parsep=\smallskipamount %
\itemsep=0.0pt \topsep=0.0pt \partopsep=\smallskipamount}}{\end{list}}
\newcounter{iii}
\newenvironment{iii}{\begin{list}{\rm \roman{iii}) }%
{\usecounter{iii} \leftmargin=0.0pt \labelsep=0.0pt %
\listparindent=0.0pt \labelwidth=0.0pt \parsep=\smallskipamount%
 \itemsep=0.0pt \topsep=0.0pt \partopsep=\smallskipamount}}{\end{list}}

\def\rightend#1#2{{%
 \leavevmode\nobreak\hskip .5em plus 1fil
 \penalty600 \hskip 0pt plus -1filll
 \vadjust{}\nobreak\hskip 0pt plus 1filll%
 #1\parfillskip=#2\relax \par}}

\def\eop{\ifmmode\rule[-22pt]{0pt}{1pt}\ifinner\tag*{$\square$}\else\eqno{\square}\fi\else\rightend{$\square$}{0pt}\fi}

\author{
Andreas-Stephan Elsenhans${}^*$ and J\"org Jahnel${}^\ddagger$}

\date{}

\title{Kummer~surfaces and the computation of \\the Picard~group}

\begin{document}

\maketitle

\begin{abstract}
We test R.\,van Luijk's method for computing the Picard group of a
$K3$~surface.
The~examples considered are the resolutions of Kummer quartics
in~$\bP^3$.
Using~the theory of abelian varieties, in this case, the Picard group may be computed~directly. Our~experiments show that the upper bounds provided by R.\,van Luijk's method are sharp when sufficiently many primes are~used. In~fact, there are a lot of primes that yield a value close to the exact~one. However, for many but not all Kummer surfaces
$V$
of Picard rank
$18$,
we have
$\smash{\rk\Pic(V_{\overline\bbF_{\!p}}) \geq 20}$
for a set of primes of density
$\geq\! \frac12$.
\end{abstract}

\footnotetext[1]{Mathematisches Institut, Universit\"at Bayreuth, Univ'stra\ss e 30, D-95440 Bayreuth, Germany,\\
{\tt Stephan.Elsenhans@uni-bayreuth.de}, Website:\! {\tt http://www.staff\!.\!uni-bayreuth.de/$\sim$btm216}\smallskip}

\footnotetext[3]{\mbox{D\'epartement Mathematik, Universit\"at Siegen, Walter-Flex-Str.~3, D-57068 Siegen, Germany,} \\
{\tt jahnel@mathematik.uni-siegen.de}, Website: {\tt http://www.uni-math.gwdg.de/jahnel}\smallskip}

\footnotetext[1]{The first author was supported in part by the Deutsche Forschungsgemeinschaft (DFG) through a funded research~project.\smallskip}

\footnotetext[1]{${}^\ddagger$The computer part of this work was executed on the servers of the chair for Computer Algebra at the University of Bayreuth. Both authors are grateful to Prof.~M.~Stoll for the permission to use these machines as well as to the system administrators for their~support.}

\section{Introduction}

\begin{ttt}
For~a general
$K3$~surface~$V$,
the methods to compute the geometric Picard group are limited up to~now. As~shown, for example in~\cite{vL2}, \cite{Kl}, \cite{EJ1}, \cite{EJ3}, or~\cite{EJ5}, it is possible to construct
$K3$~surfaces
of rank two or four with a prescribed Picard~group. But~when a
$K3$~surface
is given, say, by an equation with rational coefficients, it is not entirely clear whether its geometric Picard rank may be determined using the methods presently~known.
\end{ttt}

\begin{ttt}
To~be concrete, one may always establish a lower bound by specifying divisors explicitly and verifying that their intersection matrix is nondegenerate. On~the other hand, for upper bounds, the method of R.\,van~Luijk is available, which is based on reduction
modulo~$p$.
It~is not at all clear whether the upper bounds provided by van Luijk's method are always~sharp.
\end{ttt}

\begin{rem}
Conjecturally,~the Picard rank of a
$K3$~surface
over~$\overline\bbF_{\!p}$
is always~even. In~particular, if, for
$V$
a
$K3$~surface
over~$\bbQ$,
$\smash{\rk \Pic(V_{\overline\bbQ})}$
is odd then there is no
prime~$p$
of good reduction such
that~$\smash{\rk \Pic(V_{\overline\bbF_{\!p}}) = \rk \Pic(V_{\overline\bbQ})}$.
Even~more, the rank
over~$\smash{\overline\bbQ}$
being even or odd, there is no obvious reason why there should exist a prime
number~$p$
such that
$\smash{\rk \Pic(V_{\overline\bbF_{\!p}})}$
is at least close
to~$\smash{\rk \Pic(V_{\overline\bbQ})}$.
\end{rem}

\begin{defi}
Let~$V$
be a
$K3$~surface
over~$\bbQ$
and
$p$
be a prime of good~reduction. Then,~we will call
$p$~{\em good\/}
if the geometric Picard rank of the reduction
modulo~$p$
does not exceed the Picard rank
over~$\overline\bbQ$
by more than~one.
\end{defi}

\begin{ttt}
In~this article, we will report on our experiments concerning van Luijk's method on a sample of Kummer~surfaces. Kummer~surfaces are particular
$K3$~surfaces
allowing a two-to-one covering by an abelian~surface. The~geometric Picard group of a Kummer surface is closely related to the N\'eron-Severi~group of the abelian~surface. In~practice, it may be computed this~way.

Nevertheless,~for testing van Luijk's method, Kummer~surfaces have big advantages. Knowing~the Picard ranks anyway, the usual question whether the lower bound or the upper bound needs to be improved, does not~appear. Further,~using the particularities of a Kummer surface, one may massively optimize the point counting~step. In~fact, it is very well possible to compute
$\rk \Pic(V_{\overline\bbF_{\!p}})$
for
primes~$p$
up
to~$10\,000$.
\end{ttt}

\begin{ttt}
Our~sample consists of the resolutions of 9452 Kummer quartics with small~coefficients. For~each of these surfaces, we computed the upper bounds that were found using the
primes~$p \leq 997$.
It~turned out that good primes existed in every~example. The~upper bounds found turned out to be equal to the geometric Picard ranks in all~cases.
\end{ttt}

\begin{qu}
Do there exist good primes for every
$K3$~surface
over~$\bbQ$?
\end{qu}

\begin{lui}
\label{Luijk}
The~geometric Picard group of a
$K3$~surface
over~$\bbQ$
is isomorphic
to~$\bbZ^n$
where
$n$~may
range from
$1$
to~$20$.
An~upper bound for the geometric Picard~rank may be computed as~follows. One~has the~inequality
$$\rk \Pic (V_{\overline\bbQ}) \leq \rk \Pic (V_{\overline\bbF_p})$$
that is true for every smooth variety
$V$
over~$\bbQ$
and every
prime~$p$
of good~reduction. This~is worked out in detail in~\cite[Remark 2.6.3]{vL1}, the main input being~\cite[Example~20.3.6]{Fu}.

Further,~for a
$K3$~surface~$\calV$
over the finite
field~$\bbF_{\!p}$,
one has the first Chern class~homomorphism
$$c_1\colon \Pic (\calV_{\overline\bbF_{\!p}}) \longrightarrow H^2_\et (\smash{\calV_{\overline\bbF_{\!p}}} \!, \bbQ_l (1))$$
into
$l$-adic~cohomology.
There~is a natural operation of Frobenius on
$H^2_\et (\smash{\calV_{\overline\bbF_{\!p}}} \!, \bbQ_l (1))$.
All~eigenvalues are of absolute
value~$1$.
The~Frobenius operation on the Picard group is compatible with the operation on~cohomology.

Every~divisor is defined over a finite extension of the ground~field. Consequently,~on the subspace
$\smash{\Pic (\calV_{\overline\bbF_{\!p}}) \!\otimes_\bbZ\! \overline\bbQ_l \hookrightarrow H^2_\et (\calV_{\overline\bbF_{\!p}} \!, \bbQ_l (1))}$,
all eigenvalues are roots of~unity. These~correspond to eigenvalues of the Frobenius on
$\smash{H^2_\et (\calV_{\overline\bbF_{\!p}} \!, \bbQ_l)}$
being of the form
$p \zeta$
for
$\zeta$
a root of~unity. One~may therefore estimate the rank of the Picard
group~$\smash{\Pic (\calV_{\overline\bbF_{\!p}})}$
from above by counting how many eigenvalues are of this particular~form.

Doing~this for one prime, one obtains an upper bound
for~$\rk \Pic(V_{\overline{\bbF}_{\!p}})$
that is always~even. The~Tate conjecture asserts that this bound is actually~sharp. For~this reason, one tries to combine information from two~primes. The~assumption that the surface would have
Picard~rank~$2r$
over~$\overline{\bbQ}$
and~$\overline{\bbF}_{\!p}$
implied that the discriminants of both Picard groups,
$\Pic(V_{\overline{\bbQ}})$
and~$\Pic(V_{\overline{\bbF}_{\!p}})$,
were in the same square~class. Note~here that reduction
modulo~$p$
respects the intersection~product. When~combining information from two primes, it may happen that one finds the rank
bound~$2r$
twice, but the square classes of the discriminants are not the~same. Then,~these data are incompatible
with Picard
rank~$2r$
over~$\overline{\bbQ}$.
There~is a rank bound
of~$(2r-1)$.
\end{lui}

\begin{rem}
There~are refinements of the method of van Luijk described in~\cite{EJ3} and~\cite{EJ5}. We~will not test these refinements~here.
\end{rem}

\begin{exa}
Let~$V$
be a
$K3$~surface
of geometric
Picard~rank~$1$.
We~denote~by
$$V^n := \bigtimes_{i=1}^n V$$
the
\mbox{$n$-fold}
cartesian~product. Then,~the Picard rank
of~$V^n$
is equal
to~$n$.

Indeed,~as we
have~$H^1(V(\bbC),\bbZ) = 0$,
the K\"unneth formula shows that
$H^2(V(\bbC)^n,\bbZ) \cong H^2(V(\bbC),\bbZ)^n$.
There~is an analogous isomorphism for cohomology with complex~coefficients, which is compatible with Hodge~structures.
I.e.,~$H^{1,1}(V(\bbC)^n) \cong H^{1,1}(V(\bbC))^n$.
The~assertion now follows from the Lefschetz theorem on
\mbox{$(1,1)$-classes}~\cite[p.\,163]{GH}.

Assuming~the Tate conjecture, one sees that the Picard rank of the reduction
of~$V^n$
at an arbitrary prime is at
least~$2n$.
This~shows that there is no good~prime. Not~knowing the decomposition
of~$V^n$
into a direct product, we could not determine its Picard~rank.
\end{exa}

\begin{conv}
Let~$V$
be a projective variety over a
field~$k$.
In~this article, unless stated otherwise, the Picard rank
of~$V$
shall always mean the geometric Picard rank, i.e., the rank
of~$\smash{\Pic(V_{\overline{k}})}$.
\end{conv}

\paragraph{The analytic discriminant -- The Artin-Tate formula.}

For~the final step in~\ref{Luijk}, one needs to know the discriminant of the Picard~lattice. One~possibility to compute this is to use the Artin-Tate~formula.

\begin{cono}[{\rm Artin-Tate}{}]
\label{AT}
Let\/
$V$
be a\/
$K3$~surface
over a finite
field\/~$\bbF_{\!q}$.
Denote~by\/
$\rho$
the rank and by\/
$\Delta$
the discriminant of the Picard group
of\/~$V\!\!$,
defined
over\/~$\bbF_{\!q}$.
Then,
$$|\Delta| = \frac{\lim\limits_{T \rightarrow q}\! \frac{\Phi(T)}{(T-q)^{\rho}}}{q^{21 -\rho} \#\!\Br(V)} \, .$$
Here,
$\Phi$
denotes the characteristic polynomial
of\/~$\Frob$
on\/~$\smash{H_\et^2(V_{\overline{\bbF}_{\!q}}\!, \bbQ_l)}$.
Finally,
$\Br(V)$
is the Brauer group
of\/~$V$.
\end{cono}

\begin{rems}
\begin{iii}
\item
The Artin-Tate conjecture was first formulated, more generally than just for
$K3$~surfaces,
as Conjecture~(C) in~\cite[p.\,426]{Ta1}.
\item
Conjecture~\ref{AT} is proven for most
$K3$~surfaces.
Most~notably, the Tate conjecture implies the Artin-Tate~conjecture~\cite[Theorem~6.1]{Mi1}. In~these cases,
$\#\!\Br(V)$
is a perfect~square.

On~its part, the Tate conjecture is proven for
$K3$~surfaces
under various additional~assumptions. For~example, it is true for elliptic
$K3$~surfaces~\cite{ASD}. For~ordinary
$K3$~surfaces,
it is known, too~\cite{NO}, but we will not need this~fact.
\item
The~Artin-Tate formula allows to compute the square class of the discriminant of the Picard group over a finite~field. No~knowledge of explicit generators is~necessary.
\end{iii}
\end{rems}

\section{Singular quartics}

Singular quartic surfaces were extensively studied by the geometers
of the 19th~century, particularly by A.\,Cayley and~E.\,E.~Kummer. For~example, the concept of a trope is due to this~period~\cite{Je}.

\begin{defi}
Let
$Q \subset \bP^3$
be any quartic~surface. Then,~by a {\em trope\/}
on~$Q$,
we mean a
plane~$E$
such that
$Q \cap E$
is a double~conic. This~is equivalent to the condition that the equation
defining~$Q$
becomes a perfect square
on~$E$.
\end{defi}

\begin{rem}
A~trope yields a singular point on the surface
$Q^\vee \subset (\bP^3)^\vee$
dual
to~$Q$.
\end{rem}

\begin{lemo}[{\rm Kummer}{}]
A~quartic surface without singular curves may have at most 16 singular~points.
\eop
\end{lemo}

\paragraph{A classical family.}
A~classification of the singular quartic surfaces with at least eight singularities of
type~$A_1$
was given by K.\,Rohn~\cite{Ro}, cf.~\cite[Chapter\,I]{Je}. In~this article, we will deal with one of the most important classical~families.

\begin{lemo}[{\rm Kummer}{}]
A three-dimensional family of quartics
in\/~$\bP^3$
such that the generic member has exactly 16 singularities of
type\/~$A_1$
and no others is given by the~equation
$$16 k x y z w - \phi^2 = 0 \, .$$
Here,
\begin{eqnarray*}
 k &:=& a^2 + b^2 + c^2 - 1 - 2abc \, , \\
 \phi &:=& x^2 + y^2 + z^2 + w^2 + 2a(yz+xw) + 2b(xz+yw) + 2c(xy+zw)
\end{eqnarray*}
for parameters\/
$a$,
$b$,
and\/~$c$.
\end{lemo}

\begin{rems}
\label{Kum}
\begin{iii}
\item
E.\,E.~Kummer introduces this family in section~10 of his report~\cite{Ku}.
\item
We~will write
$Q_{[a,b,c]}$
for the quartic corresponding to the
triple~$[a,b,c]$.

Up~to isomorphism, this surface is independent of the order
of~$a,b,c$.
Further,~there is the isomorphism
$\smash{Q_{[a,b,c]} \stackrel{\cong}{\longrightarrow} Q_{[-a,-b,c]}}$
given by
$(x:y:z:w) \mapsto ((-x):(-y):z:w)$.
\item
When~one of the coefficients is equal
to~$\pm 1$,
$Q_{[a,b,c]}$~contains
a singular~line. For~example, the surfaces
for~$a = \pm 1$
contain the singular line, given
by~$x + a w = y + a z = 0$.
\item
On~the generic fiber, there are twelve obvious singularities defined over quadratic extensions
of~$\bbQ(a,b,c)$.
These~are given by
$x=y=0, z^2 + w^2 + 2czw = 0$
and the analogous conditions with the roles of the variables~interchanged. Further,~there~are four singular points forming a Galois~orbit.
\item
On a Kummer quartic, there are 16 tropes. Four~of them are~obvious. They~are explicitly given by the coordinate~planes. Each~trope passes through six of the 16~singular~points and each singular point is contained in six tropes~\cite[Chapter\,I]{Hu}.

On~an obvious trope, the conic is of discriminant
$2abc + 1 - a^2 - b^2 - c^2 = -k$.
Thus,~these conics are nondegenerate except for the case that
$Q$
is non-reduced~itself.
\item
For~a generic Kummer quartic, every singular point
on~$Q^\vee$
comes from a~trope.
\end{iii}
\end{rems}

\section{The desingularization}

\begin{lem}
\label{des}
Let\/~$\smash{\pi\colon \widetilde{Q} \to Q}$
be the desingularization of a normal quartic
surface\/~$Q$
such that all singularities are of
type\/~$A_1$.
Then,~$\smash{\widetilde{Q}}$
is a\/
$K3$~surface.\smallskip

\noindent
{\bf Proof.}
{\em
On~the smooth part
of~$Q$,
the adjunction formula~\cite[Sec.~1.1, Example~3]{GH} may be applied as~usual. As,~for the canonical sheaf, one
has~$K_{\bP^3} = \calO(-4)$,
this shows that the invertible sheaf
$\Omega^2_{Q^\reg}$
is~trivial.
Consequently,~$K_{\widetilde{Q}}$
is given by a linear combination of the exceptional~curves.

However,~for an exceptional
curve~$E$,
we have
$E^2 = -2$.
Hence, according to the adjunction formula,
$K_{\widetilde{Q}}E = 0$,
which shows that
$K_{\widetilde{Q}}$
is~trivial. The~classification of algebraic surfaces~\cite{Be} assures that
$\smash{\widetilde{Q}}$
is either a
$K3$~surface
or an abelian~surface.

Further,~a standard application of the theorem on formal functions implies
$R^1\pi_* \calO_{\widetilde{Q}} = 0$.
Hence,~$\smash{\chi_\alg(\widetilde{Q}) = \chi_\alg(Q) = 2}$.
This~shows that
$\smash{\widetilde{Q}}$
is actually a
$K3$~surface.
}
\eop
\end{lem}

\begin{rems}
\begin{iii}
\item
For~the assertion of the lemma, it is actually sufficient to assume that the singularities
of~$Q$
are of types
$A$,
$D$,
or~$E$~\cite{Li}.
\item
In~general, the desingularization of a normal quartic surface is a
$K3$~surface,
a rational surface, a ruled surface over an elliptic curve, or a ruled surface over a curve of genus~three~\cite{IN}. The~last possibility is caused by a quadruple~point. The~existence of a triple point implies that surface is~rational. It~is, however, also possible that there is a double point, not of type
$A$,
$D$,
or~$E$.
Then,~$\smash{\widetilde{Q}}$
is rational or a ruled surface over an elliptic~curve.
\end{iii}
\end{rems}

\begin{lem}
Let\/~$\smash{\pi\colon \widetilde{Q} \to Q}$
be the desingularization of a proper
surface\/~$Q$
having only\/
\mbox{$A_1$-singularities}.

\begin{abc}
\item
Then,~the exceptional curves define a nondegenerate orthogonal system
in\/~$\smash{\Pic(\widetilde{Q})}$.
\item
In~particular, the Picard~rank
of\/~$\smash{\widetilde{Q}}$
is strictly bigger than the number of singularities
of\/~$Q$.
\end{abc}\smallskip

\noindent
{\bf Proof.}
{\em
a)
The~exceptional curves have self-intersection
number~$(-2)$
and do not meet each~other.\smallskip

\noindent
b)
For~$H$
the hyperplane section,
$\pi^* \calO_Q(H)$
is orthogonal to the exceptional~curves.
}
\eop
\end{lem}

\section{Abelian surfaces and Kummer quartics}

Let~$A$ be an abelian~surface.
Denote~by
$\phi \colon A \rightarrow A$
the involution given
by~$p \mapsto (-p)$.
Then,~the
quotient~$A / \!\!\sim$
for~$\sim\,\,\, := \{(p, \phi(p)) \mid p \in A \}$
has precisely 16 singular~points. We~call such a quotient an {\em abstract Kummer~surface.}

\begin{fac}
Let\/~$A$
be an abelian surface over a
field\/~$k$
of characteristic zero
and\/
$V$~be
the resolution of the corresponding Kummer~surface.
Then,~$\smash{\rk\Pic(V_{\overline{k}}) = \rk\NS(A_{\overline{k}}) + 16}$.\smallskip

\noindent
{\bf Proof.}
{\em
A~standard argument~\cite[Proposition~(8.9.1)]{EGAIV} allows us to assume that
$k$~is
finitely generated
over~$\bbQ$.
Then,~in particular,
$k$
allows an embedding
into~$\bbC$.
The~canonical injection
$\iota\colon H^2(V(\bbC), \bbZ) \to H^2(A(\bbC), \bbZ)$
yields a bijection
of~$H^2(V(\bbC), \bbZ)$
with~$\langle E_1, \ldots, E_{16} \rangle^\perp$.
As~$\iota$
respects the
$(1,1)$-classes,
the assertion~follows. Observe~that base~change
to~$\bbC$
does not change the Picard and N\'eron-Severi~ranks.
}
\eop
\end{fac}

\begin{lemo}[{\rm Nikulin}{}]
Let\/~$Q$
be a quartic surface over an algebraically closed
field\/~$k$
of characteristic zero with precisely 16~singular
points of
type\/~$A_1$
and no~others.
Then,~$Q$
is isomorphic to an abstract Kummer~surface.\smallskip

\noindent
{\bf Proof.}
{\em
This~result is shown in~\cite{Ni}. We~include a sketch of the proof for the convenience of the~reader.

Again,~we may assume that
$k$
is a subfield
of~$\bbC$.
As~shown in Lemma~\ref{des}, the
desingularization~$\smash{\widetilde{Q}}$
is a
$K3$~surface.
We~have to prove that
$\smash{\widetilde{Q}}$
admits a double cover ramified exactly at the 16 exceptional
curves~$E_1, \ldots, E_{16}$.
This~is equivalent to the assertion that
$\smash{\calO(E_1 + \cdots + E_{16}) \in \Pic(\widetilde{Q})}$
is divisible by~two.

Consider,~more generally, the
set~$C$
of all
\mbox{$\bbQ$-divisors}
$D = c_1 E_1 + \cdots + c_{16} E_{16}$
that define an element
of~$\smash{\Pic(\widetilde{Q})}$.
Clearly,~$c_1, \ldots, c_{16} \in \frac12 \bbZ$
as, otherwise, the intersection numbers with
$E_1, \ldots, E_{16}$
would not be~integers.
Thus,~$C$
defines
a sub-vector space
$\overline{C}$~of
$$\textstyle{\bigoplus\limits_{i=1}^{16} \frac12 \bbZ E_i \,\big/ \bigoplus\limits_{i=1}^{16} \bbZ E_i \,\cong\, \bbF_{\!2}^{16} \, .}$$
We~claim that
$\dim \overline{C} \geq 5$.
Indeed,~otherwise, the lattice
$\smash{C \subset \Pic(\widetilde{Q})}$
would have a basis containing twelve of the standard elements
$E_1, \ldots, E_{16}$.
As~the quotients
$\smash{H^2(\widetilde{Q}(\bbC), \bbZ) / \Pic(\widetilde{Q})}$
and~$\smash{\Pic(\widetilde{Q}) / C}$
have no torsion,
$\smash{H^2(\widetilde{Q}(\bbC), \bbZ)}$
still has a basis containing twelve of
the~$E_i$.
But~then, the
\mbox{$22 \times 22$-matrix}
of the cup product form contains a symmetric
\mbox{$12 \times 12$-block}
consisting entirely of even~entries. This~ensures that the determinant is even and, hence, is a contradiction to the unimodularity
of~$\smash{H^2(\widetilde{Q}(\bbC), \bbZ)}$.

Further,~every vector
in~$\overline{C}$
is a sum of precisely eight or 16 standard basis~vectors. In~fact, if it is a
sum of
$l$
basis~vectors then it defines a double
cover~$P^\prime$
of~$\smash{\widetilde{Q}}$
ramified at exactly
$l$
of the 16 exceptional
curves~$E_1, \ldots, E_{16}$.
Its~minimal
model~$P$,
obtained by blowing down the
$l$
exceptional curves, clearly has trivial canonical~class. It~is therefore either an abelian surface,
$\chi_\top(P) = 0$,
or
a~$K3$~surface,
$\chi_\top(P) = 24$.
But~a direct calculation shows
$\chi_\top(P) = 48 - 3l$.

Finally,~it is a well-known result from coding theory~\cite[Theorem~2.7.4]{HP} that there is no five-dimensional subspace
of~$\bbF_{\!2}^{16}$
such that every non-zero vector has exactly eight components equal
to~$1$.
Indeed,~adding the vector
$(1,1,\ldots,1)$
would yield a code contradicting the optimality of the
$[16,5,8]$-Hadamard~code.
}
\eop
\end{lemo}

\begin{ttt}
Consider~the particular case that
$A = J(C)$
is the Jacobian of a
curve~$C$
of genus~two. Then,~a projective model of the corresponding Kummer~surface may be obtained as~follows.\smallskip

\noindent
For~$r$
a Weierstra{\ss} point
of~$C$,
put
$\theta := \{[x] - [r] \mid x \in C \} \subset J(C)$.
This~is an ample divisor on the
Jacobian~$J(C)$
such
that~$\theta^2 = 2$.
The~Riemann-Roch theorem
shows~$\dim \Gamma(J(C), 2\theta) = 4$.
Hence,~$2\theta$
defines a morphism
$\iota\colon J(C) \to \bP^3$
of degree~eight.
Actually,~$\iota$
is a two-to-one~map inducing an embedding
of
$J(C) / \!\!\!\sim$
\cite[Chapter\,VIII, Exercise\,4]{Be}.
The~image
of~$\iota$
is a quartic~surface.
\end{ttt}

\begin{ttt}
\label{g2const}
It~is a classical result that every Kummer
quartic~$Q$
may be constructed from a
\mbox{genus-$2$}~curve~$C$
in this~way. We~may therefore ask for an explicit construction of such a curve from a given Kummer~quartic. This~may indeed be done as~follows.\medskip

\begin{iii}
\item[{\bf Construction.} i) ]
There are 16~tropes. We~choose one of them, which we
call~$D$.
\item[ii) ]
The intersection
$Q \cap D$
is a double~conic.
Let~$I$
be the underlying reduced~curve. Six~of the singular points
on~$Q$
are contained
in~$I$.
\item[iii) ]
Take~the double
cover~$C$
of~$I$
ramified at these six~points. This~is a \mbox{genus-$2$~curve}.
\end{iii}
\end{ttt}

\begin{rems}
\begin{iii}
\item
This~construction clearly yields a
\mbox{genus-$2$~curve}
$C$
on the abelian
surface~$A$.
The~Albanese property of the Jacobian guarantees that
$A$~is
at least isogenous
to~$J(C)$.
They~are actually isomorphic to each~other.
\item
If~$Q$
is defined over a base field
$k$
and
$D$ over an extension
$k^\prime \supseteq k$
then
$C$
is defined
over~$k^\prime$.
Indeed,~the six ramification points form a
$\Gal(\overline{k^\prime}/k^\prime)$-invariant~set.
We~will apply the construction only to the obvious tropes of the Kummer family, which are defined over the base~field.
\end{iii}
\end{rems}

\begin{fac}
Let\/~$V^\prime$
be an abstract Kummer surface over a finite
field\/~$\bbF_{\!q}$
and\/
$V$
its resolution of~singularities. Then,~the\/
$\Gal(\overline\bbF_{\!q}/\bbF_{\!q})$-module\/
$\smash{H_\et^2(V_{\overline{\bbF}_{\!q}}\!, \bbQ_l)}$
is~reducible. A~direct summand is isomorphic
to\/~$\smash{H_\et^2(A_{\overline{\bbF}_{\!q}}\!, \bbQ_l)}$
for~$A$
the abelian~surface
covering\/~$V^\prime$.
Its~complement is described by the Galois operation on the 16 singular~points.
\eop
\end{fac}

\begin{remo}[{\rm Frobenius eigenvalues for Kummer surfaces}{}]
In~order to determine the eigenvalues of the Frobenius
on~$\smash{H_\et^2(V_{\overline{\bbF}_{\!q}}\!, \bbQ_l)}$,
the usual method is to count the points
on~$V$
defined
over~$\bbF_{\!q}$
and some of its extensions and to apply the Lefschetz trace~formula~\cite[Chapter\,VI, Theorem~12.3]{Mi2}.

For~Kummer surfaces, there is, however, a far better~method. In~fact, 16 eigenvalues are determined by the operation of Frobenius on the 16 singular~points.
Further,~for~$A$
isogenous to the
Jacobian~$J(C)$,
we have
$H_\et^2(A_{\overline{\bbF}_{\!q}}\!, \bbQ_l) \cong \Lambda^2 H_\et^1(C_{\overline{\bbF}_{\!q}}\!, \bbQ_l)$.
Thus,~in order to determine the remaining six eigenvalues, it suffices to count the points
on~$C$.
This~is faster as the problem is reduced to dimension~one.
\end{remo}

\begin{prop}
Let\/~$A$
be an abelian surface over an algebraically closed~field. Suppose~that\/
$\End(A)$
is an order of a real quadratic number~field.
Then,~$\rk \NS(A) = 2$.\smallskip

\noindent
{\bf Proof.}
{\em
According~to~\cite[sec.\,21, Appl.\,III]{Mu1}, one has
$\NS(A) \!\otimes\! \bbQ \cong (\End(A) \!\otimes\! \bbQ)^\dagger$
where
$\dagger$
denotes the Rosati~involution. As~that is positive~\cite[sec.\,21, Theorem\,1]{Mu1}, it cannot be the conjugation on a real quadratic number~field.
Hence,~$\dagger = \id$
which implies the~assertion.
}
\eop
\end{prop}

\begin{rem}
If~the real multiplication is by an order
in~$\bbQ(\sqrt{d})$
then the discriminant of the N\'eron-Severi lattice is of square
class~$(-d)$.
Indeed,~\cite[sec.\,21, Thm.\,1]{Mu1} shows that, for
$H$~ample,
$\Phi^*\!(H) \!\cdot\! H$
is a scalar multiple
of~$\Tr(\Phi^2)$.
Working~with
$\Phi = 1$
and~$\Phi = 1 + \sqrt{d}$,
we find the intersection~matrix
$\smash{\big(\atop{\,\,\,\,\,\,2 \,\,\,\,\,\,\,\,2(1+d)}{2(1+d) \,2(1-d)^2}\!\big)}$
of
determinant~$(-16d)$.
\end{rem}

\begin{prop}
\label{real}
Let\/~$A$
be an abelian surface
over\/~$\bbQ$.
Suppose~that\/
$A$
has an
endomorphism\/~$N$
defined only over a quadratic
extension\/~$F = \bbQ(\sqrt{D})$.
Then,~for every
prime\/~$p$,
inert
in\/~$F$
and of good reduction, the following is~true.\smallskip

\noindent
If\/~$\lambda$
is an eigenvalue
of\/~$\Frob_p$
on\/~$\smash{H^1_\et (A_{\overline\bbF_{\!p}}, \bbQ_l)}$
then\/
$(-\lambda)$
is an eigenvalue,~too.\medskip

\noindent
{\bf Proof.}
{\em
$N$~induces
an endomorphism
of~$A_{\bbF_{\!p}}$,
which we denote
by~$\underline{N}$.
Clearly,~$\underline{N}$
is defined
over~$\bbF_{\!p^2}$
but not
over~$\bbF_{\!p}$.
This~means, in the endomorphism
ring~$R_p$
of~$\smash{A_{\overline\bbF_{\!p}}}$,
we have
$\Frob_{p^2}^{-1} \underline{N} \Frob_{p^2} = \underline{N}$
but the analogous statement
is not true
for~$\Frob_p$.

Thus,~under the operation
of~$\Frob_{p^2}$
on~$R_p$
by conjugation,
$\underline{N}$
lies in the
\mbox{$(+1)$-eigenspace}.
For~the corresponding operation
of~$\Frob_{p}$,
this space decomposes into a
\mbox{$(+1)$-eigenspace}
and
a~\mbox{$(-1)$-eigenspace}.
The~latter is nonzero as
$\underline{N}$
is not fixed under conjugation
by~$\Frob_p$.
Hence,~there is some
$J \in R_p$
anticommuting
with~$\Frob_p$.
This~implies the~assertion.
}
\eop
\end{prop}

\begin{coro}
\label{anti}
Let\/~$V$
be a Kummer surface
over\/~$\bbQ$
covered by the abelian
surface\/~$A$.
Suppose~that\/
$A$
has an
endomorphism\/~$N$
defined only over a quadratic
extension\/~$F = \bbQ(\sqrt{D})$.\smallskip

\noindent
Then,~for every
prime\/~$p$,
inert
in\/~$F$
and of good reduction,
$\rk \Pic(V_{\overline\bbF_{\!p}}) \geq 20$.\medskip

\noindent
{\bf Proof.}
{\em
Recall~that Kummer surfaces are elliptic~\cite[Chapter\,IX, Exercise\,6]{Be}. Hence,~the Tate conjecture is true
for~$\smash{V_{\overline\bbF_{\!p}}}$.

Under~the assumptions made, it is possible that the Frobenius eigenvalues
on~$\smash{H^1_\et (A_{\overline\bbF_{\!p}}, \bbQ_l)}$
are~$\smash{\pm \sqrt{p\mathstrut}}$
and~$\smash{\pm i\sqrt{p\mathstrut}}$.
This~yields Picard
rank~$22$
over~$\overline\bbF_{\!p}$.
Except~for this case, the Frobenius eigenvalues must
be~$\pm \lambda$
and~$\pm \overline\lambda$
for a
suitable~$\lambda \in \bbC$.
On~$\smash{H^2_\et (A_{\overline\bbF_{\!p}}, \bbQ_l)}$,
this leads to the eigenvalues
$p$
and~$(-p)$,
both with multiplicity two, as well as
$(-\lambda^2)$
and~$(-\overline\lambda {}^2)$.
The~Picard rank is at
least~$20$.
}
\eop
\end{coro}

\section{The tetrahedroid}

The~tetrahedroid is another family of quartic surfaces studied in the 19th~century. It~was first considered by A.\,Cayley in~\cite{Ca1}.

\begin{lem}
A~family of quartics
in\/~$\bP^3$
such that every member has exactly 16
\mbox{$A_1$-singularities}
and no others is given by the~equation
$$
\det
\left(
\begin{array}{ccccc}
    0 &    x_0^2 &    x_1^2 &    x_2^2 &    x_3^2 \\
x_0^2 &        0 & a_{01}^2 & a_{02}^2 & a_{03}^2 \\
x_1^2 & a_{01}^2 &        0 & a_{12}^2 & a_{13}^2 \\
x_2^2 & a_{02}^2 & a_{12}^2 &        0 & a_{23}^2 \\
x_3^2 & a_{03}^2 & a_{13}^2 & a_{23}^2 &   0
\end{array}
\right) = 0 \, .
$$
for parameters\/
$a_{01}, a_{02}, a_{03}, a_{12}, a_{13}, a_{23} \neq 0$.
\end{lem}

\begin{rems}
\begin{iii}
\item
In~this form, the equation of the tetrahedroid appears in~\cite[p.~286]{Ca2}.
\item
We~will write
$T_{[a_{01}, a_{02}, a_{03}, a_{12}, a_{13}, a_{23}]}$
for the quartic corresponding to the particular coefficient
vector~$(a_{01}, a_{02}, a_{03}, a_{12}, a_{13}, a_{23})$.
\item
Let~the group
$\bbG_m^4$
operate on the parameters according to the~rule
$$(i,j,k,l)[a_{01}, a_{02}, a_{03}, a_{12}, a_{13}, a_{23}] := [ija_{01}, ika_{02}, ila_{03}, jka_{12}, jla_{13}, kla_{23}] \, .$$
Then,~the quartics defined by a whole orbit are all isomorphic to each other, as one can see by left and right multiplying the matrix above
by~$\diag(1, i^2, j^2, k^2, l^2)$.
Consequently,~the tetrahedroid defines only a two-dimensional family in the moduli stack of all
$K3$~surfaces.
Actually,~it is a subfamily of the Kummer quartics~\cite{Ca2}, \cite[\S56]{Hu}.
\end{iii}
\end{rems}

\begin{rems}
\begin{abc}
\item
The sixteen singularities are
$(0          : \pm a_{01} : \pm a_{02} : \pm a_{03})$,
$(\pm a_{01} :          0 : \pm a_{12} : \pm a_{13})$,
$(\pm a_{02} : \pm a_{12} :          0 : \pm a_{23})$,
and
$(\pm a_{03} : \pm a_{13} : \pm a_{23} :          0)$.
\item
The~four planes, given by
$\pm a_{23} x_1 \pm a_{13} x_2 \pm a_{12} x_3 = 0$,
clearly contain six singular points~each. For~example,
$((\pm a_{01}) :           0 :    a_{12} : (-a_{13}))$,
$((\pm a_{02}) :      a_{12} :         0 : (-a_{23}))$,
and
$((\pm a_{03}) :      a_{13} :  (-a_{23}) :        0)$
satisfy the
equation~$a_{23} x_1 + a_{13} x_2 + a_{12} x_3 = 0$.
There~are twelve more tropes obtained in an analogous manner by distinguishing the first, second, or third coordinate instead of the zeroth~one.
\item
Besides the tropes, there are four other particular planes related to this family of~quartics. Actually,~the coordinate planes contain exactly four singularities~each. As~these form a tetrahedron, they gave this family its~name. There~are no planes containing exactly four singular points on a general Kummer~quartic.
\end{abc}
\end{rems}

\begin{prop}
\label{tet_ell}
Let\/~$E_1$
and\/~$E_2$
be two elliptic curves. Fix~an isomorphism of groups\/
$\phi\colon E_1[2] \to E_2[2]$
and let
$$A := (E_1 \times E_2) / \langle (x, \phi(x)) \mid x \in E_1[2] \rangle$$
be the corresponding abelian surface, covered four-to-one
by\/~$E_1 \times E_2$.\smallskip

\noindent
Then,~the Kummer surface corresponding
to\/~$A$
is given by a tetrahedroid.\medskip

\noindent
{\bf Proof.}
{\em
We~describe the elliptic curves as intersections of two quadrics
in~$\bP^3$,
\begin{eqnarray*}
E_1\colon & \,x_1^2 = x_0^2 - x_2^2, \quad x_3^2 = x_0^2 - \kappa_1 x_2^2 \, , \\
E_2\colon & \,y_1^2 = y_0^2 - y_2^2, \quad y_3^2 = y_0^2 - \kappa_2 y_2^2 \, .
\end{eqnarray*}
We~have
$j(E_1) = 256 (\kappa_1^2 - \kappa_1 + 1)^3 / \kappa_1^2 (\kappa_1-1)^2$
and the analogous formula
for~$E_2$.
Thus,~these equations define general families of elliptic~curves. The~morphism
$$E_1 \times E_2 \longrightarrow \bP^3, \quad ((x_0 : x_1 : x_2 : x_3),(y_0 : y_1 : y_2 : y_3)) \mapsto (x_2y_3 : x_1y_1 : x_3y_2 : x_0y_0)$$
is generically eight-to-one onto the
tetrahedroid~$T_{[\sqrt{\kappa_2 - 1}, 1, \sqrt{\kappa_2}, \sqrt{\kappa_1 - 1}, i, \sqrt{\kappa_1}]}$.
It~factors through
$A$
and even through the Kummer surface associated with~it.
}
\eop
\end{prop}

\begin{rems}
\begin{iii}
\item
It is not hard to see that every tetrahedroid is obtained from two elliptic curves in this~way.
\item
From~the point of view of the present article, Proposition~\ref{tet_ell} is a purely algebraic~statement. We~even checked the assertions on the morphism using~{\tt magma}. It~was, however, originally discovered by H.\,Weber~\cite[p.\,353]{We} in the guise of a parametrization of the tetrahedroid by elliptic~functions. Cf.~\cite[Chapter\,XVIII]{Hu}.
\end{iii}
\end{rems}

\begin{propo}[{\rm Kummer quartics with two coefficients equal}{}]
\label{Kumeq}
\leavevmode

\noindent
Let\/~$V := V_{[a,a,c]}$
be the Kummer quartic for the
coefficients\/~$[a,a,c]$.
Then,~$V$
is linearly isomorphic to the~tetrahedroid
$$T_{\left[\sqrt{c+1}, \sqrt{c-1}, X\sqrt{c-1}, X\sqrt{c-1}, \sqrt{c-1}, 2X(X+a)\sqrt{c-1}\right]} \, .$$
Here,~$X$
is a solution of the
equation\/~$X^2 + 2aX + 1 = 0$.\smallskip

\noindent
{\bf Proof.}
{\em
The~isomorphism from the tetrahedroid
to~$V$
is given explicitly by the linear map
$\bP^3 \to \bP^3$,
$$\textstyle
(t_1 : t_2 : t_3 : t_4) \mapsto
\big(
(-t_2 - \frac{t_3-Xt_4}{1-X^2}) :
(-t_1 + \frac{t_4-Xt_3}{1-X^2}) :
(t_2 - \frac{t_3-Xt_4}{1-X^2}) :
(t_1 + \frac{t_4-Xt_3}{1-X^2})
\big) \, . \eop$$
}%
\end{propo}

\begin{rem}\vskip-3mm
\label{jinv}
One~might ask to determine the two elliptic
curves~$E_1, E_2$
which correspond
to~$V_{[a,a,c]}$,
i.e.,~those
satisfying~$(E_1 \times E_2) / \langle (x, \phi(x)) \rangle \cong V_{[a,a,c]}$.
This~leads to a simple calculation 
but the explicit formulas become rather~lengthy. Interestingly,~the two
\mbox{$j$-invariants}
are defined in the quadratic field extension
$\bbQ(a,c)(\!\sqrt{4a^2 - 2c -2})$
and conjugate to each~other. Their~trace~is
\begin{eqnarray*}
&&
(1024a^{10}c^2 + 2048a^{10}c + 1024a^{10} - 512a^8c^3 - 4608a^8c^2 - 7680a^8c - 3584a^8 \\[-1mm]
&&
{} + 32a^6c^4 + 1568a^6c^3 + 7776a^6c^2 + 10976a^6c + 4736a^6 - 72a^4c^4 - 1680a^4c^3 \\[-1mm]
&&
{} - 6016a^4c^2 - 7280a^4c - 2872a^4 + 54a^2c^4 + 702a^2c^3 + 2010a^2c^2 + 2106a^2c \\[-1mm]
&&
{} \textstyle + 760a^2 - \smash{\frac{27}2}\!c^4 - 81c^3 - 180c^2 - 175c - \smash{\frac{125}2}) / (a-1) (a+1) (b-1)^2 (b+1)^2 ,
\end{eqnarray*}
\noindent
while their norm turns out to~be
$$
\frac{(16 a^4b^2 + 48 a^4 - 24 a^2b^2 - 32 a^2b - 72 a^2 + 9 b^2 + 30 b + 25)^3}{16 (a-1)^2 (a+1)^2 (b-1)^4 (b+1)^2} \, .$$
\end{rem}

\begin{rems}
\begin{iii}
\item
The~case of three equal coefficients is even more~special. In~some sense, the
quartics~$V_{[a,a,a]}$
are tetrahedroids in three distinct~ways.

It~turns out that, in this situation, the resulting elliptic curves are related by an isogeny of
order~$3$.
In~fact, it is easy to check that the resulting pair of
\mbox{$j$-invariants}
is a zero of the third classical modular~polynomial.

Consequently,~the Picard rank of a Kummer surface with three equal coefficients is at
least~$19$.
The~additional divisor leading to a Picard rank higher
than~$18$
is the image of the graph of the
\mbox{$3$-isogeny}
under the two-to-one covering described in Proposition~\ref{tet_ell}.
\item
There~is another case, which is particular. Consider~the
quartics~$V_{[0,0,c]}$.
Then,~the 
\mbox{$j$-invariants}
of the corresponding elliptic curves are defined
in~$\bbQ(c)$
and equal to each~other. We~have
$\smash{j(E_1) = j(E_2) =  \frac{1728c^3 + 8640c^2 + 14400c + 8000}{c^3 - c^2 - c + 1}}$.

Consequently,~the Picard rank of a Kummer surface with two coefficients zero is at
least~$19$.
The~additional divisor leading to a Picard rank higher
than~$18$
is the image of the~diagonal.
\end{iii}
\end{rems}

{\boldmath
\section{Experiments -- The Picard ranks
over~$\overline\bbQ$}
}

\paragraph{A sample of Kummer~surfaces.}

%
\looseness-1
We~inspected the
Kummer~surfaces~$Q_{[a,b,c]}$
given by the Kummer~coefficients
$a, b, c = -30, \ldots, 30$.
Because~of symmetry, the considerations were restricted to the
case~$|a| \leq b \leq c$.
Recall~that one may always change the signs of two coefficients~simultaneously. Hence,~$b,c \geq 0$
was~assumed. The~coefficient vectors
$[ 3, 3, 17 ]$,
$[ 2, 2, 7 ]$,
and~$[ 2, 7, 26 ]$
as well as those containing
$\pm1$
were excluded from the sample as the corresponding surfaces have singularities of types worse
than~$A_1$.

\begin{ttt}
For~each
surface~$Q$
in the sample, first, using Construction~\ref{g2const}, we determined the
genus-$2$~curve~$C$
such that
$V$
is the Kummer surface corresponding
to~$J(C)$.
Then,~for every prime number
below~$1000$,
we counted the numbers of points
on~$C$
over
$\bbF_{\!p}$
and~$\bbF_{\!p^2}$.
From~these data, we computed the characteristic polynomial of the Frobenius on the
\mbox{$l$-adic}
cohomology of the
resolution~$V$.

From~the characteristic polynomial, we read off the rank
of~$\Pic(V_{\overline\bbF_{\!p}})$
and, using the Artin-Tate~formula~\ref{AT}, computed the square class of the~discriminant. Note~that the Artin-Tate formula is~applicable, since every Kummer surface is~elliptic.
\end{ttt}

\paragraph{The Picard ranks over~\boldmath$\overline\bbQ$.}

A~generic Kummer surface is of geometric Picard
rank~$17$.
In~the case that two Kummer~coefficients are of the same absolute value, Proposition~\ref{Kumeq} shows together with Remark~\ref{Kum}.ii) that the surface is a~tetrahedroid. Then,~the Picard~rank is at
least~$18$.
Thus,~we distinguished between these two~cases. The~possibilities that all three coefficients coincide, at least up to symmetry, or that two coefficients vanish were treated as being somehow~exceptional.

Being~a bit sloppy at first, in the first case, we tested whether an upper bound
of~$17$
is provable by van Luijk's method while, in the second case, we awaited an upper bound
of~$18$.
The~table below shows the distribution of the biggest prime that had to be considered in order to prove the~expectation.

\begin{table}[H]
\scriptsize
\begin{multicols}{2}
\hfill
\columnsep0pt
\begin{tabular}{|c|c|c|}
\hline
prime & \#cases finished & \#cases left \\
\hline\hline
\phantom{0}7  & \phantom{00}57  &           7656 \\\hline
11            & \phantom{0}287  &           7369 \\\hline
13            & \phantom{0}713  &           6656 \\\hline
17            &           1229  &           5427 \\\hline
19            &           1308  &           4119 \\\hline
23            &           1215  &           2904 \\\hline
29            &           1004  &           1900 \\\hline
31            & \phantom{0}759  &           1141 \\\hline
37            & \phantom{0}551  & \phantom{0}590 \\\hline
41            & \phantom{0}320  & \phantom{0}270 \\\hline
43            & \phantom{0}143  & \phantom{0}127 \\\hline
47            & \phantom{00}59  & \phantom{00}68 \\\hline
53            & \phantom{00}28  & \phantom{00}40 \\\hline
59            & \phantom{00}17  & \phantom{00}23 \\\hline
61            & \phantom{000}6  & \phantom{00}17 \\\hline
67            & \phantom{000}3  & \phantom{00}14 \\\hline
73            & \phantom{000}1  & \phantom{00}13 \\\hline
83            & \phantom{000}1  & \phantom{00}12 \\
\hline
\end{tabular}
\columnbreak

%
\begin{tabular}{|c|c|c|}
\hline
prime  & \#cases finished & \#cases left \\
\hline\hline
\phantom{00}5 &           156  &           1495 \\\hline
\phantom{00}7 & \phantom{0}66  &           1429 \\\hline
\phantom{0}11 &           193  &           1236 \\\hline
\phantom{0}13 &           253  & \phantom{0}983 \\\hline
\phantom{0}17 &           288  & \phantom{0}695 \\\hline
\phantom{0}19 &           132  & \phantom{0}563 \\\hline
\phantom{0}23 &           117  & \phantom{0}446 \\\hline
\phantom{0}29 &           116  & \phantom{0}330 \\\hline
\phantom{0}31 & \phantom{0}82  & \phantom{0}248 \\\hline
\phantom{0}37 & \phantom{0}81  & \phantom{0}167 \\\hline
\phantom{0}41 & \phantom{0}73  & \phantom{00}94 \\\hline
\phantom{0}43 & \phantom{0}24  & \phantom{00}70 \\\hline
\phantom{0}47 & \phantom{0}18  & \phantom{00}52 \\\hline
\phantom{0}53 & \phantom{0}15  & \phantom{00}37 \\\hline
\phantom{0}59 & \phantom{0}13  & \phantom{00}24 \\\hline
\phantom{0}61 & \phantom{00}6  & \phantom{00}18 \\\hline
\phantom{0}67 & \phantom{00}3  & \phantom{00}15 \\\hline
\phantom{0}71 & \phantom{00}2  & \phantom{00}13 \\\hline
\phantom{0}73 & \phantom{00}4  & \phantom{000}9 \\\hline
\phantom{0}79 & \phantom{00}2  & \phantom{000}7 \\\hline
          101 & \phantom{00}1  & \phantom{000}6 \\
\hline
\end{tabular}
\hfill
\end{multicols}\vskip-4mm
\caption{Distribution of the biggest prime used for rank 17 (left) and rank 18 (right)}
\normalsize
\end{table}

\paragraph{The 18 examples left.}

Let~us take a closer look at the Kummer quartics~left.

\begin{exas}
Among~the Kummer quartics, the coefficients of which had three distinct absolute values, twelve examples~remained. For~these, only a rank bound
of~$18$
could be~established. Using~{\tt magma}, we calculated the corresponding
\mbox{genus-$2$~curves}~$C_i$
and determined their periods at high~precision.

\begin{iii}
\item
Consider~the Kummer quartics for the coefficient
vectors
$[ 2, 3, 13 ]$,
$[ -3, 4, 19 ]$,
$[ -3, 5, 11 ]$,
$[ -2, 7, 23 ]$,
$[ -2, 8, 17 ]$,
$[ -2, 9, 14 ]$,
and~$[ 0, 4, 7 ]$.

In~these cases, it turned out that the Jacobians
$J(C_i)$
are isogenous to products of two elliptic~curves. Hence, the geometric Picard ranks are indeed equal
to~$18$.

The~isogenies are all of
degree~$16$.
Their~kernels are groups of
type~$\bbZ/4\bbZ \times \bbZ/4\bbZ$.
The~\mbox{$j$-invariants}
of the elliptic curves are conjugate to each other in quadratic number~fields. We~summarize them in the table~below.

\begin{table}[H]
\scriptsize
\begin{center}
\begin{tabular}{|c|r@{$\pm$}l|}
\hline
vector  &  \multicolumn{2}{c|}{$j_1, j_2$} \\
\hline\hline
\rule{0pt}{2.4ex}
$[ \phantom{-}2, \phantom{1}3,           13 ]$ & $\frac{8000}{21609}[38155$ & $16152 \sqrt{2}]$           \\[0.5mm]\hline
\rule{0pt}{2.4ex}
$[           -3, \phantom{1}4,           19 ]$ & $\frac{64}{164025}[1082783$ & $399784 \sqrt{-2}]$        \\[0.5mm]\hline
\rule{0pt}{2.4ex}
$[           -3, \phantom{1}5,           11 ]$ & $\frac{16}{5625}[17903$ & $64596 \sqrt{-1}]$             \\[0.5mm]\hline
\rule{0pt}{2.4ex}
$[           -2, \phantom{1}7,           23 ]$ & $\frac{2000}{10673289}[-614135$ & $4744012 \sqrt{-2}]$   \\[0.5mm]\hline
\rule{0pt}{2.4ex}
$[           -2, \phantom{1}8,           17 ]$ & $\frac{250}{21609}[-50045$ & $45683 \sqrt{-3}]$          \\[0.5mm]\hline
\rule{0pt}{2.4ex}
$[           -2, \phantom{1}9,           14 ]$ & $\frac{16}{6426225}[327552721$ & $229629540 \sqrt{-1}]$  \\[0.5mm]\hline
\rule{0pt}{2.4ex}
$[ \phantom{-}0, \phantom{1}4, \phantom{1}7 ]$ & $\frac{16}{5625}[17903$ & $64596 \sqrt{-1}]$             \\[0.5mm]
\hline
\end{tabular}
\end{center}\vskip-2mm
\caption{$j$-invariants of the corresponding elliptic curves}
\normalsize
\end{table}

\item
Consider~the Kummer quartics given by the coefficient vectors
$[ 2, 7, 17 ]$,
$[ 2, 9, 26 ]$,
$[ 2, 17, 26 ]$,
$[ 3, 9, 19 ]$,
and~$[ 0, 8, 15 ]$.

Here,~our calculations showed that the corresponding abelian surfaces have real multiplication by orders
in~$\bbQ(\sqrt{2})$,
$\bbQ(\sqrt{3})$,
$\bbQ(\sqrt{5})$,
$\bbQ(\sqrt{5})$,
and~$\bbQ(\sqrt{5})$,~respectively.
This~implies that the Picard ranks are equal
to~$18$.

The~non-trivial endomorphisms are expected to be defined over the quadratic number fields
$\bbQ(\sqrt{30})$,
$\bbQ(\sqrt{11})$,
$\bbQ(\sqrt{-2})$,
$\bbQ(\sqrt{-1})$,
and~$\bbQ(\sqrt{2})$.
In~fact, the primes leading to Picard
rank~$18$
are all split for the corresponding~field. Compare~Corollary~\ref{anti}.
\end{iii}
\end{exas}

\begin{ttt}
Consider~the Kummer quartics for the coefficient vectors
$[ 5, 5, 17 ]$,
$[ 2, 2, 17 ]$,
$[ -4, 4, 9 ]$,
$[ -3, 7, 7 ]$,
$[ -2, 11, 11 ]$,
and~$[ 0, 5, 5 ]$.
For~these, the situation is as~follows.\smallskip

\noindent
One~finds 
rank~$20$
at several~primes. Discriminants~of various square classes appear such that a rank bound
of~$19$
is~established.\smallskip

\noindent
As~two of the Kummer coefficients are equal, the corresponding abelian surfaces are isogenous to products of two elliptic~curves. Specializing~the calculation discussed in Remark~\ref{jinv}, one may determine the corresponding
\mbox{$j$-invariants.}
It~turns out in every case that the corresponding elliptic curves are isogenous to each~other. Thus,~we have Picard
rank~$19$.
The~isogenies are of degrees
$5$,
$5$,
$4$,
$4$,
$4$,
and
$4$,~respectively.
\end{ttt}

\begin{exas}
Let~us present two of these examples in~detail.

\begin{iii}
\item
Let~$T_1$
be the Kummer quartic for the coefficient vector
$[5, 5, 17]$.
We~find
rank~$20$
at
$p = 5$,
$7$,
$13$,
$17$,
$19$,
$23$,
$29$
and several other~primes.
The~rank
bound~$19$
is proven as many distinct square classes of discriminants~occur.

Further,~as two of the Kummer coefficients are equal, the corresponding abelian surface is isogenous to a product of two elliptic~curves. Specializing~the calculation discussed in Remark~\ref{jinv}, one finds the
\mbox{$j$-invariants}
$j_1 = \frac{85184}{3}$
and~$j_2 = \frac{58591911104}{243}$.
The~pair
$(j_1, j_2)$
is a zero of the fifth modular~polynomial. Hence,~between the two elliptic curves, there is an isogeny of order~five. We~have Picard
rank~$19$.
\item
Let~$T_2$
be the Kummer quartic for the coefficient vector
$[2, 2, 17]$.
For~this surface, we find
rank~$20$
at
$p = 7$,
$11$,
$13$,
$17$,
$23$,
$29$
and several other~primes. Many~distinct square classes of discriminants~appear.
Hence,~the rank
bound~$19$
is~proven.

Here,~the two
\mbox{$j$-invariants}
are defined
in~$\bbQ(\sqrt{-5})$.
They~are the roots of the
polynomial~$X^2 + \frac{21180800}{243}X + \frac{1693669888000}{729}$.
Again,~the corresponding elliptic curves turn out to be
\mbox{$5$-isogenous}.
This~confirms Picard
rank~$19$.
\end{iii}
\end{exas}

\paragraph{Expected rank 19.}

In~the case of three equal coefficients or two coefficients equal to zero,
we know that the Picard rank is at least equal
to~$19$.
In~84 of the 88 surfaces, the reductions
modulo~$p$
provided an upper bound
of~$19$.
The~biggest prime that had to be used
was~$37$.

The~cases
$[0, 0, 0]$,
$[-5, 5, 5]$,
$[-2, 2, 2]$,
and
$[7, 7, 7]$~remained.
Here,~the corresponding elliptic curves have complex~multiplication. This~shows that the corresponding Kummer surfaces indeed have geometric Picard
rank~$20$.

\begin{exa}
Consider,~for instance, the case
$[7, 7, 7]$.
Then,~the two
\mbox{$j$-invariants}
are the roots of the polynomial
$X^2 - 37018076625X + 153173312762625$.
The~corresponding elliptic curves have complex multiplication by an order
in~$\bbQ(\sqrt{-15})$.
\end{exa}

\paragraph{Testing isomorphy.}

As~a byproduct of the computations, we tried to prove that the surfaces in our sample are pairwise non-isomorphic. For~this, it would suffice to show that,
for each pair of surfaces, there exists a prime where both have good reduction, but the geometric Picard groups differ in rank or~discriminant. Actually,~the data
for~$p \leq 59$
contained enough information for this but there were 41~pairs of surfaces that could not be~separated.

The point here is that the test actually tries to prove that the corresponding abelian surfaces are non-isogenous. But~in these 41 cases, the surfaces
are isogenous to each~other. To~be more precise, we found 17 pairs, four triples, and two quadruples of mutually isogenous abelian~surfaces.

\begin{exa}
The~abelian surfaces corresponding to
$V_{[2,2,9]}$
and~$V_{[3,3,19]}$
are~isogenous. Hence,~the test described above has no chance to~work.

In~fact,
$V_{[2,2,9]}$
is covered eight-to-one
by~$E_1 \times E_2$
while
$V_{[3,3,19]}$
is covered eight-to-one
by~$E_3 \times E_4$
for
$j(E_1), j(E_2)$
the zeroes
of~$X^2 + \frac{1114112}{25}X + \frac{589752696832}{225}$
and
$j(E_3), j(E_4)$
the zeroes of
$\smash{X^2 - \frac{281615072}{2025}X + \frac{15000601854041872}{164025}}$.
It~is easy to check that
$E_1$
and~$E_3$,
as well as
$E_2$
and~$E_4$,
are connected by isogenies of order~four.
Hence,~$E_1 \times E_2$
and~$E_3 \times E_4$
are~\mbox{$16$-isogenous}.

An~isomorphism between the quotients as described
in~Proposition~\ref{tet_ell} would yield
a~\mbox{$16$-isogeny}
\begin{eqnarray*}
& & \hspace{-0.2cm} E_1 \times E_2 \longrightarrow E_1 \times E_2 / \langle (x, \phi(x)) \mid x \in E_1[2] \rangle \cong E_3 \times E_4 / \langle (x, \phi^\prime(x)) \mid x \in E_3[2] \rangle \\
& & \hspace{12cm} \longrightarrow E_3 \times E_4 \, ,
\end{eqnarray*}
too. But,~in its kernel, there are the
\mbox{$2$-torsion}
points
$(x, \phi(x))$
for~$x \in E_1[2]$
that are clearly not in the kernel of the direct product of two
\mbox{$4$-isogenies}.
This~shows that
$V_{[2,2,9]}$
and~$V_{[3,3,19]}$
are not isomorphic,~either.
\end{exa}

\paragraph{Testing isomorphy II.}
For~each of the 41~pairs, we numerically calculated the periods of the corresponding abelian~surfaces. From~these, we determined a minimal~isogeny. It~turned out that the surfaces corresponding to the coefficient vectors
$[ -3, 7, 7 ]$
and~$[ 0, 5, 5 ]$
were actually isomorphic to each~other. This~was, however, the only such case among the critical~pairs.

\paragraph{Summary.}
We~considered the resolutions of 9452 Kummer~quartics with exactly 16 singularities of
type~$A_1$.
It~turned out that the upper bounds for the Picard ranks provided by the reductions
modulo~$p$
were sharp in every~case. However,~at several examples, rather large primes up
to~$p = 101$
had to be~considered. We~had
Picard
rank~$17$,
$7701$~times,
Picard
rank~$18$,
$1657$~times,
and Picard
rank~$19$,
$90$~times.
Further,~there were four
surfaces of Picard
rank~$20$
in the~sample.

\section{Some more statistics}

\begin{exo}[{\rm All primes less than
$10\,000$
for a typical surface}{}]
\leavevmode

\noindent
Let~us take a closer look at a particular~example. We~selected the surface with
Kummer
coefficients~$[3, 11, 21]$,
but many others would be representative, as~well.

There~are only five primes
$p \leq 10\,000$
such that the reduction
modulo~$p$
of~$V_{[3,11,21]}$
is not a quartic having 16 singular points of
type~$A_1$.
These~are
$2$,
$3$,
$5$,
$11$,
and~$17$.
In~the range considered, 1224~primes lead to a reduction of Picard
rank~$18$.
Further,~there are 69 primes leading to a reduction of
rank~$20$.
These~seem to be rather equidistributed within the range, the smallest one
being~$7$,
the largest one
being~$9677$.
Finally,~there is the prime
$4583$
that leads to a reduction of Picard
rank~$22$.

In~the cases of reduction to
Picard rank~$18$,
we found 586 different square classes for the~discriminant. As~for many of the surfaces in our sample, the most frequent square class
was~$(-1)$.
In~the example selected, it appeared
$376$~times.
\end{exo}

\paragraph{Discriminants -- The special case of rank 17.}

In~the special case of a
rank-$17$~surface,
we counted how many square classes of discriminants occurred when reducing to surfaces of Picard
rank~$18$
modulo various~primes. There~are
$168$~prime
numbers in our computational range, e.g.~less
than~$1000$.
For~a fixed surface, between
$44$
and~$89$
distinct square classes were~found.

\begin{figure}[H]
\centerline{
\includegraphics[height=5.0cm]{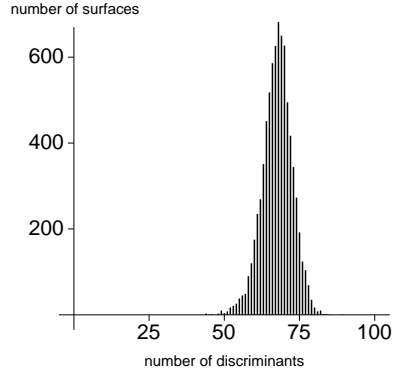}}\vskip-3mm
\caption{Number of distinct square classes of discriminants at primes with reduction to rank 18}
\end{figure}\vskip-2mm

\noindent
In~total, we found 541 distinct square classes of~discriminants. Some~of them occurred only for one surface and one~prime. On~the other hand, the class
of~$(-1)$
appeared
$134\,553$~times.
The~surfaces with Kummer coefficients
$[ -3, 9, 17 ]$
and~$[ -3, 10, 29 ]$
both had the most repetitions for one square~class. This~was the class
of~$(-1)$
occurring
$43$~times.

\paragraph{The average value for a prime.}

For~simplicity, let us restrict our considerations to surfaces of Picard
rank~$17$.
For~every prime
number~$p$,
we counted how many of the surfaces in our sample had good reduction
modulo~$p$.
We~determined the proportion of those having reduction to
rank~$>\!\!18$.
The~results are visualized by the graph~below. According to this graph, the proportion is close
to~$\frac{C}{\sqrt{p}}$
for~$C \approx 2$.

\begin{figure}[H]
\centerline{
\includegraphics[height=5cm]{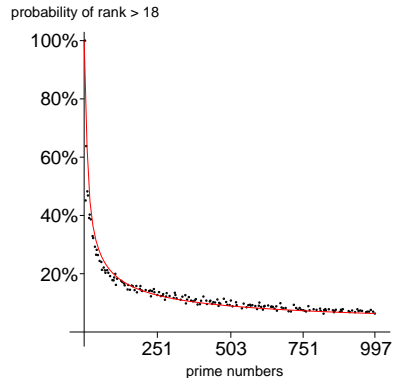}}\vskip-3mm
\caption{Distribution of the proportion of the surfaces with reduction to rank $>\!18$}
\label{Zwei}
\end{figure}

\paragraph{The average value for a surface.}

On~the other hand, for every surface of Picard
rank~$\leq \!18$
in the sample, we counted how many primes
below~$1000$
lead to a reduction of geometric Picard
rank~$18$
over~$\overline{\bbF}_{\!p}$.
Let~us visualize the result in a~histogram.

\begin{figure}[H]
\centerline{
\includegraphics[height=5cm]{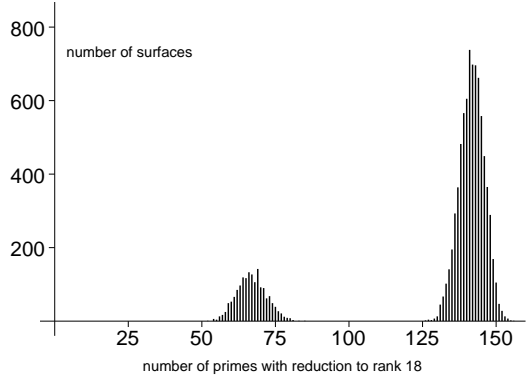}}\vskip-3mm
\caption{Distribution of the number of the primes with reduction to rank 18}
\end{figure}

\noindent
The~histogram clearly suggests that there are two kinds of~examples.
For~the first kind, the probability that the reduction has
rank~$18$
is between
$1/4$
and~$1/2$.
For~the second kind, this probability is between
$3/4$
and~$1$.

It~turns out that most of the examples with two Kummer coefficients equal (up to sign) belong to the first~kind. The~only examples in the first group not being of this form are given by the coefficient vectors
$[ 3, 9, 19 ]$,
$[ 2, 3, 13 ]$,
$[ 2, 7, 17 ]$,
$[ 2, 9, 26 ]$,
$[ 2, 17, 26 ]$,
$[ -3, 4, 19 ]$,
$[ -3, 5, 11 ]$,
$[ -2, 7, 23 ]$,
$[ -2, 8, 17 ]$,
$[ -2, 9, 14 ]$,
$[ 0, 4, 7 ]$,
and~$[ 0, 8, 15 ]$.
Further,~there are some examples with two coefficients equal belonging to the second~group. These~are
$[ 3, 3, 9 ]$,
$[ 3, 3, 15 ]$,
$[ 4, 4, 13 ]$,
$[ 4, 4, 23 ]$,
$[ 4, 4, 29 ]$,
$[ 6, 6, 21 ]$,
$[ 7, 7, 25 ]$,
$[ 8, 8, 29 ]$,
$[ 2, 2, 5 ]$,
$[ -6, 6, 27 ]$,
$[ -5, 5, 23 ]$,
$[ -4, 4, 19 ]$,
$[ -3, 3, 15 ]$,
$[ -2, 2, 11 ]$,
and~$[ -2, 2, 25 ]$.

\paragraph{An explanation.}

For~the tetrahedroid case, an explanation is given by the following~fact.

\begin{fac}
\label{tetconj}
Let~$V_{[a,a,c]}$
be a Kummer surface with two coefficients~equal.
Suppose~that
$4a^2 - 2c -2$
is not a perfect~square.\smallskip

\noindent
Then,~for every
prime\/~$p$,
inert
in\/~$F = \bbQ(\sqrt{4a^2 - 2c -2})$
and of good reduction,
$\rk \Pic(V_{\overline\bbF_{\!p}}) \geq 20$.\medskip

\noindent
{\bf Proof.}
{\em
The~corresponding abelian surface is isogenous to the product of two elliptic~curves. As~noticed in Remark~\ref{jinv}, the
\mbox{$j$-invariants}
are two elements conjugate
in~$F$.
Reducing~the surface modulo a prime inert
in~$F$
leads to two elliptic curves isogenous via the Frobenius~endomorphism. This~shows that all inert primes yield an upper bound of at
least~$20$
for the geometric Picard~rank.
}
\eop
\end{fac}

\begin{qus}
\begin{iii}
\item
For~a
surface~$V$,
put~$N_V(B) := \#\{p \in P_V \mid p \leq B\}$,
where
$$P_V := \#\{\, p {\rm ~prime} \mid \rk\Pic(V_{\overline\bbF_{\!p}}) > 18 {\rm ~or~} V {\rm ~has~bad~reduction~at~} p \,\} \, .$$
Is~there a monotonically decreasing
function~$h_V$
such~that
$$N_V(B) \sim \int\limits_2^B\! \frac{h_V(t)}{\log t} \,dt \, ?$$
Can~$h_V$
be given~explicitly?
\item
Suppose~that
$\rk\Pic(V_{\overline\bbQ}) = 17$.
Does~then
$h_V$
converge
to~$0$
for~$t \to \infty$?
The~graph in~Figure~\ref{Zwei} might suggest
that~$\smash{h_V(t) = \frac{C_V}{\sqrt{t}}}$
for a
constant~$C_V$.
Is~$h_V$
perhaps independent
of~$V$?
\item
For a fixed Kummer surface of geometric Picard
rank~$17$,
are there infinitely many primes with reduction to
rank~$18$?
Are~there infinitely many primes with reduction to
rank~$> \!18$?
\end{iii}
\end{qus}

\begin{rem}
In~relation with these questions, the reader might want to consult~\cite{MP}, for example Conjecture~5.1 formulated~there.
\end{rem}

\begin{rem}
When~$\rk\Pic(V_{\overline\bbQ}) = 18$,
the situation is typically~different. For~example, when two Kummer coefficients are equal, we saw in~Fact~\ref{tetconj} that
$\smash{P_{V_{[a,a,c]}}}$
has density at
least~$\frac12$
unless~$4a^2 - 2c - 2$
is a perfect~square.
According~to Proposition~\ref{real}, the same is true when the abelian surface corresponding
to~$V$
has real multiplication by an endomorphism defined over a proper field extension
of~$\bbQ$.
Note~that the latter case actually subsumes the former as the abelian variety corresponding to
$V_{[a,a,c]}$
is isogenous to the product of two elliptic curves and, therefore, has real~multiplication.
\end{rem}

\section{Our data}

\begin{ttt}
The~raw data of our experiments are available from NSF's Data Conservancy project as the file {\tt kummer.tar.gz} associated with this~article. They~are also available on both authors' web~pages.
\end{ttt}

\end{document}